\documentclass[reqno,11pt]{amsart}

\usepackage{amssymb,amsfonts, amsmath}
\usepackage{mathrsfs}
\usepackage{color}
\usepackage[normalem]{ulem} 
\usepackage{pstricks, pst-plot, floatflt}

\newtheorem{thm}{Theorem}[section]

\newtheorem{lem}[thm]{Lemma}

\theoremstyle{definition}

\theoremstyle{remark}

\newtheorem{exmp}[thm]{Example}
\numberwithin{equation}{section}
\newcommand{\Real}{\mathbb R}
\newcommand{\Borel}{\mathscr B}
\newcommand{\eps}{\varepsilon}

\newcommand{\F}{\mathscr{F}}

\newcommand{\one}[1]{\mathbf{1}_{\{#1\}}}

\renewcommand{\P}{\mathsf{P}}

\newcommand{\E}{\mathsf{E}}
\newcommand{\simplex}{\mathcal{S}^{d-1}}

\DeclareMathOperator*{\argmax}{\mathrm{argmax}}

\newcommand{\D}{\mathscr{D}}

\newcommand{\s}{\mathbb{S}}
\newcommand{\M}{\mathcal{M}}
\newcommand{\J}{\mathcal{J}}
\definecolor{gray}{rgb}{0.9,0.9,0.9}

\input cyracc.def

\begin{document}

\title[Filter stability for slow Markov chains]{Stability of the nonlinear filter for slowly switching Markov chains}%
\author{Pavel Chigansky}
\address{Department of Mathematics, The Weizmann Institute of Science, Rehovot 76100, Israel}
\email{pavel.chigansky@weizmann.ac.il}

\thanks{Research supported by a grant from the Israel Science Foundation}%

\date{6, July, 2006}%
\dedicatory{Dedicated to Robert Liptser on the occasion of his 70th birthday}%
\begin{abstract}
Exponential stability of the nonlinear filtering equation is revisited, when
the signal is a finite state Markov chain. An asymptotic upper bound
for the filtering error due to incorrect initial condition is derived in the case
of slowly switching signal.
\end{abstract}
\maketitle
\section{Introduction and the main result }
\label{sec-1}
Consider a discrete time Markov chain $X=(X_n)_{n\in \mathbb{Z}_+}$ with values in a finite real
alphabet $\mathbb{S}=\{a_1,...,a_d\}$, initial distribution $\nu_i= \P(X_0=a_i)$ and transition
probabilities $\lambda_{ij}=\P(X_n=a_j|X_{n-1}=a_i)$. Suppose that the chain is partially observed via the noisy
sequence of random variables $Y=(Y_n)_{n\in\mathbb{Z}_+}$, generated by
\begin{equation}
\label{Y}
Y_n = \sum_{i=1}^d \one{X_n=a_i}\xi_n(i), \quad n\ge 1,
\end{equation}
where $\xi=(\xi_n)_{n\ge 1}$ is a sequence of i.i.d. random vectors with independent entries $\xi_n(i)$, $i=1,...,d$,
such that
$$
\P\big(\xi_1(i)\in B\big)=\int_B g_i(u)\varphi(du), \quad B\in\Borel(\Real)
$$
with densities $g_i(u)$ and a $\sigma$-finite reference measure $\varphi(du)$.

Let $\F^Y_n=\sigma\{Y_1,...,Y_n\}$ and $\pi_n(i)=\P(X_n=a_i|\F^Y_n)$. The vector $\pi_n$ of the conditional probabilities
satisfies the recursive {\em filtering} equation
\begin{equation}
\label{filt}
\pi_n = \frac{G(Y_n)\Lambda^*\pi_{n-1}}{\big|G(Y_n)\Lambda^*\pi_{n-1}\big|}, \quad \pi_0=\nu,
\end{equation}
where $G(y)$, $y\in\Real$ is a diagonal matrix with entries $g_i(y)$, $\Lambda^*$ is the transposed matrix of transition probabilities and
$|x|=\sum_{i=1}^d |x_i|$ for $x\in\Real^d$.

Suppose that \eqref{filt} can be solved subject to a probability distribution $\bar{\nu}\ne \nu$ and denote the corresponding solution by $\bar{\pi}_n$. Under
certain mild conditions (to be specified later) the limit
$$
\gamma :=\lim_{n\to\infty}\frac{1}{n}\log |\pi_n-\bar{\pi}_n|, \quad \P-a.s.
$$
exists and if it is negative the filter is said to be (exponentially) stable. The {\em stability index} $\gamma$
is elusive for explicit calculation and much research focused recently on estimating $\gamma$ in various filtering
settings (see  \cite{AZ97a,AZ97b,BCL,CL,DZ,DG,LGM} and others).
In particular, Gaussian additive white noise model was considered in \cite{AZ97a}
(cf. \eqref{Y})
$$
Y_n=h(X_n)+\sigma\eta_n, \quad n\ge 1, \quad \eta_1\sim \mathcal{N}(0,1)
$$
and the following asymptotic upper bound was derived
\begin{equation}
\label{AZ}
\varlimsup_{\sigma\to 0}\sigma^2\gamma(\sigma)\le -\frac{1}{2}\sum_{i=1}^d \mu_i \min_{j\ne i}\big(h(a_i)-h(a_j)\big)^2,
\end{equation}
where $\mu$ is the stationary distribution of the chain $X$, assumed to be ergodic. Recall that $X$ is ergodic if $\mu_i=\lim_{n\to\infty}\P(X_n=a_i)$, $i=1,...,d$
exist, are unique and positive, which holds iff $\Lambda^q$ has positive entries for some integer $q\ge 1$ (see e.g. \cite{Nor}).

In this note a different scaling of the problem is chosen, namely the {\em slow chain} limit of $\gamma$ is considered.
Let $X^\eps_n$ be a Markov chain on $\s$ with transition probabilities
$$
\lambda^\eps_{ij} = \P(X^\eps_n=a_j|X^\eps_{n-1}=a_i)=\begin{cases}
\eps \lambda_{ij}, & i\ne j\\
1-\eps\sum_{\ell\ne i} \lambda_{i\ell}, & i=j.
\end{cases}
$$
for an $\eps\in(0,1)$. Notice that $X^\eps$ is an ergodic chain with the same invariant distribution $\mu$ as $X$. Denote by $Y^\eps$ the
corresponding observation sequence generated by \eqref{Y}, with $X$ replaced by $X^\eps$
and let $\pi^\eps$, $\bar{\pi}^\eps$ be the solutions of \eqref{filt} subject to $\nu$, $\bar{\nu}$, with $Y$ and $\Lambda$ replaced by
$Y^\eps$ and $\Lambda^\eps$.
\begin{thm}\label{thm}
Assume that $X$ is ergodic and the noise densities $g_i(u)$
\begin{enumerate}
\renewcommand{\theenumi}{a\arabic{enumi}}
\item \label{a1} are bounded
\item \label{a2} have the same support
\item \label{a3} and
$
\int_\Real g_i(u) \log g_j(u)\varphi(du)>-\infty,
$ for all $i,j$.
\end{enumerate}
Then for any pair $(\nu, \bar{\nu})$ of probability distributions on $\s$
\begin{equation}
\label{mr}
\gamma(\eps)\le -\sum_{i=1}^d \mu_i \min_{j\ne i}\D(g_i\parallel g_j) +o(1), \quad \eps\to 0,
\end{equation}
where $\D(g_i\parallel g_j)=\int_\Real g_i(u)\log\dfrac{g_i}{g_j}(u) \varphi(du)$ are the Kullback-Leibler relative entropies. For $d=2$ the
asymptotic \eqref{mr} is precise, i.e.
\begin{equation}\label{2d}
\gamma(\eps) = -\mu_1 \D(g_1\parallel g_2) -\mu_2\D(g_2\parallel g_1) +o(1), \quad \eps\to 0.
\end{equation}
\end{thm}
This theorem reveals the following interesting properties of $\gamma(\eps)$ (see Figure \ref{fig1}).

\noindent
{\bf 1.} $\gamma(\eps)$ may be discontinuous at $\eps=0$
$$
\gamma(0+)=\varlimsup_{\eps\to 0}\gamma(\eps)<\gamma(0)=0,
$$
if at least one of the entropies $\D(g_i\parallel g_j)$ is strictly positive.
This means that for small $\eps>0$ the filter remains stable virtually with {\em the same stability index} as long as the chain is not ``frozen'' completely,
while the filter, corresponding to the limit chain $X^0_n\equiv X_0$, $n\ge 1$, may be unstable (e.g. when some but not all $g_i(u)$'s coincide
$\varphi$-a.s.). Such a behavior is not observed in the analogous ``slowly varying'' setting for the Kalman-Bucy filter, where the state space of the signal
is continuous.

\begin{floatingfigure}[r]{0.4\textwidth}
\centering
\pspicture(-1.0,-4.5)(4.0,1.0)
\psaxes*[linewidth=0.5pt,ticks=none,labels=none]{->}(0,0)(0,-4)(3,0.5)
\rput(3.3,0){\footnotesize $\eps$} \rput(-0.25,0){\footnotesize
$0$} \psdots*[dotstyle=o,dotsize=2.0pt](0,-1.5)
\psdots*[dotstyle=*,dotsize=2.0pt](0,0)
\pscurve[linewidth=0.5pt]{->}(0,-1.5)(1.0,-0.5)(2.4,-4)
\psline[linewidth=0.5pt,linestyle=dotted,dotsep=1.5pt](2.5,0)(2.5,-4)
\rput(2.3,-4.2){\footnotesize $-\infty$}
\psdots*[dotstyle=*,dotsize=2.0pt](2.5,0)
\rput(2.5,0.3){\footnotesize $\frac{1}{2\lambda}$}
\psline[linewidth=0.5pt,linestyle=dotted,dotsep=1.5pt](0.95,0.0)(0.95,-0.5)
\rput(1.0,0.2){\footnotesize $\eps^\star$}
\rput(-0.55,-1.5){\footnotesize $-\D_p$}
\psline[linewidth=0.5pt,linestyle=dotted,dotsep=1.5pt](0.0,-1.5)(0.55,-1.5)
\psline[linewidth=0.5pt,linestyle=dotted,dotsep=1.5pt]{<->}(0.4,-0.9)(0.4,-1.5)
\psline[linewidth=0.5pt,linestyle=dotted,dotsep=1.5pt](0.4,-1.5)(0.4,-2.0)
\rput(1.0,-2.1){\footnotesize $\sim \eps\log\frac{1}{\eps}$}
\rput(0,0.7){\footnotesize $\gamma(\eps)$}
\endpspicture

\caption{\label{fig1}$\gamma(\eps)$ for the BSC example}
\noindent \hrulefill
\end{floatingfigure}

Surprising as it may seem at first glance, this phenomenon is quite natural for signals with discrete state space
and can be explained as follows. The distance $|\pi^\eps_n-\bar{\pi}^\eps_n|$ never increases and tends to decrease
exponentially fast whenever $X^\eps_n$ resides in a state with distinct noise probability distribution.
Since the average occupation time of this ``synchronizing'' state does not depend on $\eps$, the decay remains exponential
with nonzero average rate. The ``dual'' manifestation of this phenomenon is that the filter stability improves, when the signal-to-noise
ratio is increased in the setting of \eqref{AZ} (see \cite{DZ,AZ97a}).

{\bf 2.} As demonstrated in the following example, $\gamma(\eps)$ may have a maximum at some $\eps^\star>0$ or, in other words,
stability may improve when the chain is slowed down! This provides yet another evidence against the false intuition, directly relating
stability of the filter to ergodic properties of the signal (see an extended discussion of this issue in \cite{DZ,BCL}).
The reason for such behavior stems from the delicate interplay between two stabilizing mechanisms: ergodicity of the signal and
synchronizing effect of the observations. The first dominates the second for the faster chain, and vise versa when the chain is slow.

\begin{exmp}
Consider the so called Binary Symmetric Channel (BSC) model, for which $X_n\in\{0,1\}$ is a symmetric chain with the jump
probability $\lambda$ and $Y_n= (X_n -\xi_n)^2$, where $\xi$ is an i.i.d. $\{0,1\}$ binary sequence with $\P(\xi_1=1)=p\in(0,1/2)$. Let $X^\eps$ and $Y^\eps$
denote the ``slow'' instances as defined above. In this case more can be said about the convergence in \eqref{2d} (see the proof in Section \ref{ex-proof} below), namely
\begin{equation}
\label{ex-eq}
\gamma(\eps) \ge  -\D_p +
\frac{4\lambda\big(\log(2)-h(p)\big)}
{\D_p}\eps\log\eps^{-1}\big(1+o(1)\big), \quad \eps\to 0.
\end{equation}
where $\D_p:=p\log\dfrac{p}{1-p}+(1-p)\log \dfrac{1-p}{p}$ and $h(p)=-p\log p-(1-p)\log(1-p)$.
On the other hand, $\gamma(\eps)\le \log (1-2\eps\lambda)\to -\infty$ as $\eps\to 1/(2\lambda)$ (see e.g. Theorem 2.3 in \cite{AZ97a}).
Since the second term in the expansion of $\gamma(\eps)$ in \eqref{ex-eq} is positive and by \eqref{2d} $\gamma(\eps)\to -\D_p$ as $\eps\to 0$,
one gets the qualitative behavior depicted in Figure \ref{fig1}.\qed

\end{exmp}

\section{The proof of Theorem \ref{thm}}

Hereafter the assumptions of Theorem \ref{thm} are in force and the following notations are used: probability measures on $\s$ are identified with
(column) vectors in $\simplex=\{x\in\Real: x_i\ge 0, \sum_{i=1}^d x_i=1\}$, $\mu(f):=\sum_{i=1}^df(a_i)\mu_i$ for $f:\s\mapsto\Real$ and
$\mu\in\simplex$, $\mu(A):=\mu(\one{A})$ for $A\subseteq \s$. For a random sequence $Z=(Z_n)_{n\in\mathbb{Z}}$ and $m\ge k$
the notation  $\F^Z_{[k,m]}=\sigma\{Z_k,...,Z_m\}$ is used and $\F^Z_n:=\F^Z_{[1,n]}$ for brevity. Convergence of random sequences is understood in $\P$-a.s. sense unless stated otherwise.

The proof relies on the following idea from \cite{AZ97a}. Recall that $\pi_n=\rho_n/|\rho_n|$, $n\ge 0$ where $\rho_n$ is the solution of
Zakai linear equation ($\bar{\pi}_n$ is obtained similarly)
\begin{equation}
\label{zak}
\rho_n = G(Y_n)\Lambda^* \rho_{n-1}, \quad \rho_0=\nu.
\end{equation}
Let $\rho_n\wedge \bar{\rho}_n:=\frac{1}{2}\big(\rho_n\bar{\rho}_n^*-\bar{\rho}_n\rho_n^*\big)$ denote the {\em exterior product} of $\rho_n$ and $\bar{\rho}_n$.
The elementary inequality
$$
\frac{|\rho_n\wedge \bar{\rho}_n|}{|\rho_n||\bar{\rho}_n|}\le |\pi_n-\bar{\pi}_n|\le 2 \frac{|\rho_n\wedge \bar{\rho}_n|}{|\rho_n||\bar{\rho}_n|}
$$
implies
\begin{multline}\label{gamma}
\gamma:= \lim_{n\to\infty}\frac{1}{n}\log|\pi_n-\bar{\pi}_n| =  \lim_{n\to\infty}\frac{1}{n}\log |\rho_n\wedge \bar{\rho}_n| \\ -
\lim_{n\to\infty}\frac{1}{n}\log |\rho_n|-\lim_{n\to\infty}\frac{1}{n}\log |\bar{\rho}_n|.
\end{multline}
Since $g_i(u)$'s are bounded, the limits in the right hand side exist by virtue of the Oseledec Multiplicative Ergodic Theorem (MET).
Moreover, since $\big(G(Y_n)\Lambda^*\big)_{n\ge 1}$ are matrices with nonnegative entries, the Perron-Fro\-be\-nius theorem implies
$$
\lim_{n\to\infty}\frac{1}{n}\log |\rho_n|=\lim_{n\to\infty}\frac{1}{n}\log |\bar{\rho}_n|:=\lambda_1, \quad \forall \nu, \bar{\nu}\in\simplex,
$$
where $\lambda_1$ is the top Lyapunov exponent corresponding to $\eqref{zak}$. Similarly MET implies
$
\lim_{n\to\infty}\frac{1}{n}\log |\rho_n\wedge \bar{\rho}_n| \le \lambda_1+\lambda_2
$
and thus one concludes that $\gamma\le \lambda_2-\lambda_1\le 0$, i.e. the filter stability index is controlled by the
Lyapunov spectral gap of \eqref{zak}. The reader is referred to \cite{AZ97a} for further details.

The statement of Theorem \ref{thm} follows from \eqref{gamma} and asymptotic expressions derived in Lemmas \ref{lem1} and \ref{lem2} below.

\subsection{Asymptotic expression for $\lambda_1(\eps)$}
\begin{lem}\label{lem1}
For any $\eps>0$ the Markov process $(X^\eps,\pi^\eps)$ has a unique stationary invariant measure $\M^\eps$.
The top Lyapunov exponent is given by
\begin{equation}
\label{ave}
\lambda_1(\eps) = \int_{\simplex}\sum_{i=1}^d \big(\Lambda^{\eps *} u\big)_i\int_\Real g_i(y)\log\big|G(y)\Lambda^{\eps *} u\big|\varphi(dy)\M^\eps_\pi(du),
\end{equation}
where $\M^\eps_\pi$ is the $\pi$-marginal of $\M^\eps$. For each $\mathcal{J}_j=\{a_\ell:\D(g_j\parallel g_\ell)=0\}$
\begin{equation}\label{cntr}
\lim_{\eps\to 0}\int \big(\one{x\in \J_j}-\sum_{\ell: a_\ell\in \J_j}u_\ell\big)^2\M^\eps(dx,du)=0
\end{equation}
and in particular
\begin{equation}\label{lambda1}
\lim_{\eps\to 0}\lambda_1(\eps)= \sum_{i=1}^d \mu_i \int_\Real g_i(y)\log g_i(y)\varphi(dy).
\end{equation}
\end{lem}
\begin{proof}
The process $(X^\eps,\pi^\eps)$ is Markov and by \eqref{a1} it is also Feller and
thus at least one invariant measure $\M^\eps$ exists. Its uniqueness can be deduced (as in Theorem 7.1 in \cite{BK})
from the stability property $\lim_{n\to\infty}|\pi^\eps_n-\bar{\pi}^\eps_n|=0$, $\forall \nu, \bar{\nu}\in\simplex$, which in turn holds
under the assumption \eqref{a2} by the arguments used in the proof of Theorem 2.3 in \cite{AZ97a} (see also Theorem
4.1 in \cite{BCL}). Concentration properties of $\M^\eps_\pi$ have been studied in \cite{KZ}, when all the noises are distinct, i.e.
$\D(g_i\parallel g_j)>0$ for all $i\ne j$, which is not necessarily the case here.

Let $\widetilde{X}^\eps$ be the stationary chain (i.e. $\widetilde{X}_0\sim \mu$) and $\widetilde{\pi}^\eps$  the corresponding
optimal filtering process, generated by \eqref{filt} subject to $\widetilde{\pi}^\eps_0=\mu$. For an $f:\s\to\Real$ and $n,m\ge 0$
($\widetilde{Y}^\eps$ denotes the observations corresponding to $\widetilde{X}^\eps$)\label{pagelabel1}
\begin{align*}
&\E\big(f(\widetilde{X}^\eps_{n+m})-\widetilde{\pi}^\eps_{n+m}(f)\big)^2 =
\E\Big(f(\widetilde{X}^\eps_{n+m})-\E\big(f(\widetilde{X}^\eps_{n+m})\big|\F^{\widetilde{Y}^\eps}_{n+m}\big)\Big)^2\le\\
&\E\Big(f(\widetilde{X}^\eps_{n+m})-\E\big(f(\widetilde{X}^\eps_{n+m})\big|\F^{\widetilde{Y}^\eps}_{[m+1,n+m]}\big)\Big)^2\stackrel{\dagger}{=}
\E\Big(f(\widetilde{X}^\eps_n)-\E\big(f(\widetilde{X}^\eps_n)\big|\F^{\widetilde{Y}^\eps}_n\big)\Big)^2=\\
& \E\big(f(\widetilde{X}^\eps_n)-\widetilde{\pi}^\eps_n(f)\big)^2,
\end{align*}
where stationarity of $(\widetilde{X}^\eps,\widetilde{Y}^\eps)$ have been used in $\dagger$. This means that the
filtering error for the stationary signal does not increase with time. Then
by uniqueness of $\M^\eps$ for any fixed $m\ge 0$
\begin{multline}\label{kuku}
\int \big(f(x)-u(f)\big)^2\M^\eps(dx,du)= \\
\lim_{n\to\infty}\E\big(f(\widetilde{X}^\eps_n)-\widetilde{\pi}^\eps_n(f)\big)^2\le
\E\big(f(\widetilde{X}^\eps_m)-\widetilde{\pi}^\eps_m(f)\big)^2.
\end{multline}
Define
$$
\widehat{\pi}^\eps_n(i) = \frac{\mu_i\prod_{k=1}^n g_i(\widetilde{Y}^\eps_k)}
{
\sum_{j=1}^d \mu_j\prod_{k=1}^n g_j(\widetilde{Y}^\eps_k)
}, \quad i=1,...,d
$$
and let  $A^\eps_m=\{\widetilde{X}^\eps_k=\widetilde{X}_0,\ \forall k\le m\}$, the event that $\widetilde{X}^\eps_k$ does not jump on $[0,m]$.
Notice that on the set $A^\eps_m$, the observation process is independent of $\eps$, namely
$$
\widetilde{Y}^\eps_k \equiv \widetilde{Y}^0_k=\sum_{i=1}^d \one{\widetilde{X}_0=a_i}\xi_k(i), \quad k=1,...,m.
$$
Then by optimality of $\widetilde{\pi}^\eps$
\begin{align*}
&\E\big(f(\widetilde{X}^\eps_m)-\widetilde{\pi}^\eps_m(f)\big)^2\le
\E\big(f(\widetilde{X}^\eps_m)-\widehat{\pi}^\eps_m(f)\big)^2=\\
&
\E\one{A^\eps_m}\big(f(\widetilde{X}_0)-\widehat{\pi}^0_m(f)\big)^2 +
\E\one{\Omega\setminus A^\eps_m}\big(f(\widetilde{X}^\eps_m)-\widehat{\pi}^\eps_m(f)\big)^2\le \\
&\E \big(f(\widetilde{X}_0)-\widehat{\pi}^0_m(f)\big)^2 + 4d^2 \max_{a_i\in\s}|f(a_i)|^2\big(1-\P(A^\eps_m)\big)\xrightarrow[\eps\to 0]{}\E\big(f(\widetilde{X}_0)-\widehat{\pi}^0_m(f)\big)^2
\end{align*}
For $f(x):=\one{x\in\J_j}$ the latter and \eqref{kuku} implies
$$
\varlimsup_{\eps\to 0}\int \big(\one{x\in \J_j}-\sum_{\ell: a_\ell\in \J_j}u_\ell\big)^2\M^\eps(dx,du) \le
\E\big(f(\widetilde{X}_0)-\widehat{\pi}_m(f)\big)^2\xrightarrow[m\to\infty]{}0,
$$
where the convergence holds since $\{\widetilde{X}_0\in \J_j\}\in \F^{\widetilde{Y}^0}_\infty=\bigvee_{n\ge 1}\F^{\widetilde{Y}^0}_n$
by definition of $\J_j$ and since $\widehat{\pi}^0_m(i)$, $i=1,...,d$ are the optimal estimates of $\one{\widetilde{X}_0=a_i}$ given $\F^{\widetilde{Y}^0}_m$.

Once the existence of ergodic stationary pair $(X^\eps,\pi^\eps)$ is established\footnote{such pair can be generated by taking both $X_0$ and
$\pi_0$ randomly distributed according to $\M^\eps$ and its definition can be extended to the negative times by the usual arguments.
Note that this is different from $(\widetilde{X}^\eps,\widetilde{\pi}^\eps)$ used in the proof
of $\M^\eps$ concentration }
one may use it to realize the limit $\lambda_1$  by means of the approach due to H.Furstenberg and R.Khasminskii (see e.g. \cite{Kh}). The idea is to study the
growth rate of $\rho^\eps_n$ by projecting it on the unit sphere ($\simplex$ in this case):
$$
|\rho^\eps_n| = \big|G(Y^\eps_n)\Lambda^{\eps *}\rho^\eps_{n-1} \big|=|\rho^\eps_{n-1}|\Big|G(Y^\eps_n)\Lambda^{\eps *}\frac{\rho^\eps_{n-1}}{|\rho^\eps_{n-1}|} \Big|=
|\rho^\eps_{n-1}|\big|G(Y^\eps_n)\Lambda^{\eps *}\pi^\eps_{n-1} \big|.
$$
Then by the law of large numbers (LLN) for ergodic processes (the required integrability conditions are provided by \eqref{a1} and \eqref{a3})
\begin{align}
&\lambda_1(\eps) =  \lim_{n\to\infty}\frac{1}{n}\log |\rho^\eps_n|= \lim_{n\to\infty}\frac{1}{n}\sum_{m=1}^n\log \big|G(Y^\eps_n)\Lambda^{\eps *}\pi^\eps_{n-1} \big|=
\E \log\big|G(Y^\eps_1)\Lambda^{\eps *}\pi^\eps_{0} \big| =
\nonumber
\\
&\E \sum_{i=1}^d \one{X^\eps_1=a_i}\log\big|G\big(\xi_1(i)\big)\Lambda^{\eps *}\pi^\eps_{0} \big|=
\E \sum_{i=1}^d \P\big(X^\eps_1=a_i|\F^{Y^\eps}_{(-\infty,0]}\big)\log\big|G\big(\xi_1(i)\big)\Lambda^{\eps *}\pi^\eps_{0} \big|=\nonumber\\
&\E \sum_{i=1}^d \big(\Lambda^{\eps *}\pi^\eps_0\big)_i\log \big|G\big(\xi_1(i)\big)\Lambda^{\eps *}\pi^\eps_{0} \big|.\label{onehas}
\end{align}
The latter expression is nothing but \eqref{ave}. The asymptotic \eqref{lambda1} follows from $\Lambda^\eps=I+O(\eps)$ and
the concentration \eqref{cntr} of $\M^\eps$ as $\eps\to 0$, since $g_i(u)$'s coincide $\varphi$-almost surely for all $a_i\in\J_j$ for any $j$
and the $X$-marginal of $\M^\eps$ is given by $\M^\eps_X(dx)=\sum_{i=1}^d \mu_i\delta_{a_i}(dx)$.\qed
\end{proof}

\subsection{Asymptotic bound for $\lambda_1(\eps)+\lambda_2(\eps)$}
\begin{lem}\label{lem2}  For any $\nu, \bar{\nu}\in \simplex$
\begin{multline}\label{LL}
\lim_{n\to\infty}\frac{1}{n}\log |\rho^\eps_n\wedge \bar{\rho}^\eps_n| \le \\
\sum_{i=1}^d \mu_i \max_{k\ne m}\int_\Real g_i(u)\log\big(g_m(u)g_k(u)\big)\varphi(du)+o(1),\quad\eps\to 0.
\end{multline}
In the case $d=2$
\begin{multline}\label{2dLL}
\lim_{n\to\infty}\frac{1}{n}\log |\rho^\eps_n\wedge \bar{\rho}^\eps_n| =
\log(1-\eps\lambda_{12}-\eps\lambda_{21}) +\\ \mu_1 \int_\Real g_1(u)\log\big(g_1(u)g_2(u)\big)\varphi(du)
+\mu_2 \int_\Real g_2(u)\log\big(g_1(u)g_2(u)\big)\varphi(du).
\end{multline}
\end{lem}
\begin{proof}
The process $Q^\eps_n:=\rho^\eps_n\wedge \bar{\rho}^\eps_n$ evolves in the space of antisymmetric matrices (with zero diagonal) and satisfies the
linear equation
$$
Q^\eps_n = G(Y^\eps_n) \Lambda^{\eps *}Q^\eps_{n-1}\Lambda^\eps G(Y^\eps_n), \quad Q^\eps_0 = \nu\wedge \bar{\nu},
$$
or in the componentwise notation
$$
Q^\eps_n(i,j) = \sum_{1\le k\ne \ell \le d} g_k(Y^\eps_n)\lambda^\eps_{ki}Q^\eps_{n-1}(k,\ell)\lambda^\eps_{\ell j}g_\ell(Y^\eps_n), \quad i\ne j.
$$
Unlike in the case of \eqref{zak}, it is not clear whether the limit $\lim_{n\to\infty}\frac{1}{n}\log|Q^\eps_n|$ depends on $\nu, \bar{\nu}$ or
$\Pi^\eps_n = Q^\eps_n/|Q^\eps_n|$ has any useful concentration properties as $\eps\to 0$. However the technique used in the previous
section still gives the upper bound. With a fixed integer $r\ge 1$
\begin{align*}
|Q^\eps_{n}| =& |Q^\eps_{n-r}|\Big|
\Big\{G(Y^\eps_n)\Lambda^{\eps *} ... \Big\{G(Y^\eps_{n-r+1})\Lambda^{\eps *}
\Pi^\eps_{n-r}\Lambda^\eps G(Y^\eps_{n-r+1})\Big\} ... \Lambda^\eps G(Y^\eps_n)\Big\}\Big|
\le \\
&|Q^\eps_{n-r}|\Big(\sum_{i\ne j}\big|\Pi^\eps_{n-r}(i,j)\big|\prod_{m=n-r+1}^n
g_i(Y^\eps_m)g_j(Y^\eps_m)
+ c_1(r) \eps
\Big)\le\\
&|Q^\eps_{n-r}|\Big(\max_{i\ne j}\prod_{m=n-r+1}^n
g_i(Y^\eps_m)g_j(Y^\eps_m)
+ c_1(r) \eps
\Big)
,\quad n\ge r
\end{align*}
with a constant $c_1(r)>0$, depending only on $r$ (due to assumption \eqref{a1}). By the MET the limit $\lim_{n\to\infty}\frac{1}{n}\log|Q^\eps_{n}|$
exists $\P$-a.s and hence (recall the definitions of $\widetilde{Y}^\eps$ and $A^\eps_r$ on page \pageref{pagelabel1})
\begin{align*}
&\lim_{n\to\infty}\frac{1}{n}\log|Q^\eps_{n}|=
\lim_{\ell\to\infty}\frac{1}{\ell r}\log|Q^\eps_{\ell r}|\le \\
& \le \lim_{\ell\to \infty}\frac{1}{\ell}\sum_{k=1}^\ell \frac{1}{r}\log
\Big(\max_{i\ne j}\prod_{m=kr-r+1}^{kr}
g_i(Y^\eps_m)g_j(Y^\eps_m)
+ c_1(r) \eps
\Big)\stackrel{\dagger}{=} \\
& \frac{1}{r}\E \log
\Big(\max_{i\ne j}\prod_{m=1}^r
g_i(\widetilde{Y}^\eps_m)g_j(\widetilde{Y}^\eps_m)
+ c_1(r) \eps
\Big)\le \\
& \frac{1}{r}\E \one{A^\eps_r}\log
\Big(\max_{i\ne j}\prod_{m=1}^r
g_i(\widetilde{Y}^\eps_m)g_j(\widetilde{Y}^\eps_m)
+ c_1(r) \eps
\Big)  + c_2(r)\big(1-\P_\mu(A^\eps_r)\big)\le  \\
&\frac{1}{r}\sum_{\ell=1}^d\mu_\ell \E \log
\Big(\max_{i\ne j}\prod_{m=1}^r
g_i\big(\xi_m(\ell)\big)g_j\big(\xi_m(\ell)\big)
+ c_1(r) \eps
\Big)  + c_3(r)\big(1-\P_\mu(A^\eps_r)\big)\xrightarrow{\eps\to 0} \\
&
\sum_{\ell=1}^d\mu_\ell \E\max_{i\ne j} \frac{1}{r}\log
\prod_{m=1}^r
g_i\big(\xi_m(\ell)\big)g_j\big(\xi_m(\ell)\big),
\end{align*}
where the LLN was used in $\dagger$ and $c_i(r)$ stand for $r$-dependent constants. Applying the LLN once again one gets for each $\ell$
\begin{multline*}
\frac{1}{r}\log
\prod_{m=1}^r
g_i\big(\xi_m(\ell)\big)g_j\big(\xi_m(\ell)\big) =
\frac{1}{r}\sum_{m=1}^r \log
g_i\big(\xi_m(\ell)\big)g_j\big(\xi_m(\ell)\big) \xrightarrow{r\to\infty}\\
\int_\Real g_\ell(u)\log \big(g_i(u) g_j(u)\big)\varphi(du),  \quad \P-a.s.
\end{multline*}
Since ``$\max$'' is a continuous function
$$
\max_{i\ne j} \frac{1}{r}\log
\prod_{m=1}^r
g_i\big(\xi_m(\ell)\big)g_j\big(\xi_m(\ell)\big) \xrightarrow{r\to\infty}
\max_{i\ne j}\int_\Real g_\ell(u)\log \big(g_i(u) g_j(u)\big)\varphi(du)
$$
and by the uniform integrability, provided by assumption \eqref{a3},
$$
\E\max_{i\ne j} \frac{1}{r}\log
\prod_{m=1}^r
g_i\big(\xi_m(\ell)\big)g_j\big(\xi_m(\ell)\big)\xrightarrow{r\to\infty}
\max_{i\ne j}\int_\Real g_\ell(u)\log \big(g_i(u) g_j(u)\big)\varphi(du).
$$
Putting all parts together one gets the bound \eqref{LL}. In the case $d=2$, the process $Q^\eps_n$ is one dimensional and all the
calculations can be carried out exactly, leading to the expression \eqref{2dLL}.\qed
\end{proof}

\section{Proof of \eqref{ex-eq}}\label{ex-proof}
When the observation process $Y^\eps_n$ takes values in a discrete alphabet $\s'=\{b_1,...,b_{d'}\}$, the
conditional densities (with respect to the point measure $\varphi(dy)=\sum_{i=1}^{d'}\delta_{b_i}(dy)$) are of the form
$$
g_i(y)=\sum_{j=1}^{d'} p_{ij}\one{y=b_j}, \quad \sum_{j=1}^{d'} p_{ij}=1, \ p_{ij}\ge 0,
$$
and hence by \eqref{onehas} ($\pi^\eps_{1|0}:=\Lambda^{\eps *}\pi^\eps_0$ for brevity)
\begin{multline}
\lambda_1(\eps)=\E \log\big|G(Y^\eps_1)\Lambda^{\eps *}\pi^\eps_{0} \big|=
\E \sum_{j=1}^{d'}\one{Y^\eps_1=b_j}\log\Big(\sum_{i=1}^d p_{ij}\pi^\eps_{1|0}(i)\Big)=\\
\E \sum_{j=1}^{d'} \P\big(Y^\eps_1=b_j|\F^{Y^\eps}_{(-\infty,0]}\big)\log \P\big(Y^\eps_1=b_j|\F^{Y^\eps}_{(-\infty,0]}\big)=: -\mathscr{H}(Y^\eps),
\label{lambdaH}
\end{multline}
where $\mathscr{H}(Y^\eps)$ is known as the entropy rate of the stationary process $Y^\eps=(Y^\eps_n)_{n\in\mathbb{Z}}$.

Consider now the special case, when $X^\eps$ and $Y^\eps$ take values in $\s=\{0,1\}$ and $p=\P(Y^\eps_n=i|X^\eps_n=j)$ for $i\ne j$.
The vector $\pi^\eps_n$ is one dimensional and hence
$
\P\big(Y^\eps_1=1|\F^{Y^\eps}_{(-\infty,0]}\big) = (1-p) \pi^\eps_{1|0} + p(1-\pi^\eps_{1|0}),
$
where
\begin{equation}\label{pi10}
\pi^\eps_{1|0}:= \P\big(X^\eps_1=1|\F^{Y^\eps}_{(-\infty,0]}\big) = (1-\eps\lambda_{10})\pi^\eps_0+\eps \lambda_{01}(1-\pi^\eps_0)
\end{equation}
and $\pi^\eps_0:=\P(X^\eps_0=1|\F^{Y^\eps}_{(-\infty,0]})$ are redefined for brevity.

Let $h(x):=-x\log x-(1-x)\log(1-x)$, $x\in[0,1]$ and $\ell_p(q)=(1-p) q + p(1-q)$, and define
$$
H(p,q) := h\big(\ell_p(q)\big) \quad p,q\in [0,1],
$$
where $0\log 0\equiv 0$ is understood. Since $h(x)\le \log(2)$ with equality at $x=1/2$ and
$\ell_p(1/2)=1/2$,  $H(p,q)\le \log(2)$ for all $p,q\in[0,1]$ with equality at $q=1/2$. Since $h(x)$ is a
concave function, symmetric around $x=1/2$
$$
H(p,q) = h\big((1-p) q + p(1-q)\big)\ge qh(1-p)+(1-q)h(p) = h(p), \quad p\in[0,1],
$$
with equality at $q=0$ and $q=1$. Finally for any fixed $p\in[0,1]$, $q\mapsto H(p,q)$ inherits concavity and
symmetry from $h(x)$. These properties imply the following lower bound
\begin{equation}\label{lbound}
H(p,q)
\ge h(p) + \frac{\log(2)-h(p)}{1/2}\min(q,1-q), \quad p,q\in [0,1].
\end{equation}
By Theorem 1 in \cite{KZ} for the symmetric chain $X^\eps$ with jump probability $\lambda$ and $p\ne 1/2$
\begin{multline}
\label{KZthm}
\E \min(\pi^\eps_0, 1-\pi^\eps_0)=\P\big(X^\eps_0\ne \argmax_{i}\pi^\eps_0(i)\big)=\\
\frac{\lambda}{\D_p}\eps\log\eps^{-1}\big(1+o(1)\big), \quad \eps\to 0,
\end{multline}
where $\D_p:=p\log\dfrac{p}{1-p}+(1-p)\log \dfrac{1-p}{p}$.
The expression for $\mathscr{H}(Y^\eps)$ in the case $d=2$ reads
$$
\mathscr{H}(Y^\eps) =\E H(p, \pi^\eps_{1|0})=\E H(p,\pi^\eps_0) + O(\eps),\quad \eps\to 0
$$
where the latter asymptotic follows from \eqref{pi10}, since $H(p,q)$ is differentiable in $q$.

Now \eqref{lbound} and \eqref{KZthm} imply
$$
\mathscr{H}(Y^\eps) \ge h(p) + 2\big(\log(2)-h(p)\big)
 \frac{\lambda}{\D_p}\eps\log\eps^{-1}\big(1+o(1)\big), \quad \eps\to 0,
$$
and \eqref{ex-eq} follows from \eqref{gamma}, \eqref{2dLL} and \eqref{lambdaH}.\qed
\section*{Acknowledgement}
The author is grateful to Rami Atar for elaborations around \cite{AZ97a} and the useful comments regarding
the results of this paper.

\end{document}